\documentclass{article}
\usepackage{graphicx} 
\usepackage{amsmath} 
\usepackage{amssymb} 
\usepackage{setspace}
\usepackage{cleveref}

\makeatletter
\@addtoreset{equation}{section}
\makeatother

\title{Properties of Solutions to the Full Fractional Heat Operator Equation}
\author{Lu Haipeng, Yu Mei}

\begin{document}

\maketitle

\begin{abstract}
In this paper, we consider the following indefinite fully fractional heat equation involving the master operator
. Under certain assumptions of 
the indefinite nonlinearity and its weight, 
we prove that there is no positive bounded solution, which is based on the monotonicity of the solution along the first direction that is proved by employing the method of moving planes. Besides, if the weight satisfy other conditions, we come to different conclusions according to the behavior of the nonlinearity at infinity.

To overcome the difficulties caused by the operator, 
we lead in some mathematics tools that, as we believe, will be useful in studying problems involving other fractional operators or nonlinearities.
\end{abstract}


\section{Introduction}


\quad \quad Space-time nonlocal equations arise in physical and biological contexts such as anomalous diffusion, chaotic dynamics, and biological invasions. They also model financial scenarios where waiting times between transactions correlate with subsequent price jumps \cite{13,11,10,12}. A key application is in continuous-time random walks, where the nonlinear term represents particle distributions subject to coupled random jumps and time delays. This model generalizes Brownian random walks by incorporating time nonlocality (capturing anomalously long waiting times) and space nonlocality, which allows for large jumps such as Lévy flights across distant regions.\cite{29}

The main objective of this paper is to study the following master equation with indefinite nonlinearity:
\begin{equation}
(\partial_t - \Delta)^s u(x,t) = a(x)f(u(x,t)), (x,t) \in \mathbb{R}^n \times \mathbb{R}, \label{1.1}\end{equation}
where the operator \((\partial_t-\Delta)^s,\) which was first introduced by M. Riesz in \cite{1}, is a nonlocal pseudo differential operator of order \(2s\) in space and of order \(s\) in time. Its definition by singular integral is 
\begin{equation} 
(\partial_t-\Delta)^s u(x,t):=C_{n,s} \textup{P.V.} \int_{-\infty}^t \int_{\mathbb{R}^n} \frac{u(x,t)-u(y,\tau)}{(t-\tau)^{\frac{n}{2}+1+s}} e^{-{\frac{|x-y|^2}{4(t-\tau)}}} \textup{d}y\textup{d}\tau. \label{1.2} 
\end{equation}
Here \(s \in (0,1)\), P.V. stands for Cauchy principal value and \(C_{n,s}\) is the normalization constant. From \eqref{1.2} we know that the value of \((\partial_t-\Delta)^s u(x,t)\) depends on the values of \(u\) in the whole space \(\mathbb{R}^n\) and at all time before \(t,\) which makes \((\partial_t-\Delta)^s\) nonlocal both in space and time.

Next we lead in the slowly increasing function space and the local parabolic H\"{o}lder space, which guarantee that the singular integral in \eqref{1.2} is well defined.

\textbf{Definition 1.1.}
The slowly increasing function space \(L(\mathbb{R}^n \times \mathbb{R})\) is defined as
\[L(\mathbb{R}^n \times \mathbb{R}) :=\left\{ u(x,t) \in L^1_{\textup{loc}}(\mathbb{R}^n \times \mathbb{R})| \int_{-\infty}^t \int_{\mathbb{R}^n} \frac{|u(x,\tau)|e^{-\frac{|x|^2}{4(t-\tau)}}}{1+(t-\tau)^{\frac{n}{2}+1+s}}  \textup{d}x\textup{d}\tau <\infty \right\}, \]
where \(L^1_{\textup{loc}}(\mathbb{R}^n \times \mathbb{R})\) stands for the space of local measurable functions.

\textbf{Definition 1.2.}
The parabolic H\"{o}lder space \(C^{2\alpha,\alpha}_{x,t}(\mathbb{R}^n \times \mathbb{R})\) is defined differently according to the value of \(\alpha.\) More precisely,\newline
(i)For \(0<\alpha \leq\frac{1}{2},\) we say that \(u(x,t) \in C^{2\alpha,\alpha}_{x,t}(\mathbb{R}^n \times \mathbb{R})\) if there exists a constant \(C>0\) such that for all \(x,y \in \mathbb{R}^n\) and for all \(t \in \mathbb{R},\)
\[|u(x,t)-u(y,\tau)| \leq C\left(|x-y|+|t-\tau|^\frac{1}{2}\right)^{2\alpha}. \]
(ii)For \(\frac{1}{2}<\alpha\leq1,\) we say that
\[u(x,t) \in C^{2\alpha,\alpha}_{x,t}(\mathbb{R}^n \times \mathbb{R}):=C^{1+(2\alpha-1),\alpha}_{x,t}(\mathbb{R}^n \times \mathbb{R}),\]
if \(u\) is \(\alpha-\)H\"{o}lder continuous in \(t\) uniformly with respect to \(x,\) and its gradient \(\nabla_x u\) is \((2\alpha-1)-\)H\"{o}lder continuous in \(x\) uniformly with respect to \(t\) and \((\alpha-\frac{1}{2})-\)H\"{o}lder continuous in \(t\) uniformly with respect to \(x.\) \newline
(iii)For \(\alpha>1,\) we say that \(u(x,t) \in C^{2\alpha,\alpha}_{x,t} (\mathbb{R}^n \times \mathbb{R}),\) if
\[\partial_t u, D^2_x u \in C^{2\alpha-2,\alpha-1}_{x,t},\]
where \(D^2_x u\) is the Hessian matrix of \(u.\)

Similarly, we may define the local parabolic H\"{o}lder space \(C^{2\alpha,\alpha}_{x,t,\textup{loc}}(\mathbb{R}^n \times \mathbb{R}).\)

\textbf{Definition 1.3.} We say that \(u\) is a classical entire solution of \eqref{1.1} if \(u\) satisfies \eqref{1.1} and
\[u(x,t) \in C^{2s+\epsilon,s+\epsilon}_{x,t,\textup{loc}}(\mathbb{R}^n \times \mathbb{R}) \cap L(\mathbb{R}^n \times \mathbb{R})\]
for some \(s \in (0,1)\) and \(\epsilon>0,\) which makes sure that \eqref{1.2} is well defined.

If the order \(s \to 1,\) then the operator \((\partial_t-\Delta)^s\) will tend to the regular heat operator \((\partial_t-\Delta)\)\cite{2}. If \(u\) is a function only related to space variables \(x\) or time variable \(t\), then the operator will become the fractional Laplacian \((-\Delta)^s\) or the Marchaud fractional derivative \(\partial_t^s\)\cite{3}. These operators, owing to their diverse scientific applications, are widely used in areas such as anomalous diffusion, quasi-geostrophic flows, phase transitions, image processing, plasma turbulence, and many others. They also play a significant role in probability and finance, where they can be interpreted as the infinitesimal generators of stable Lévy processes\cite{4,5,6,7,8,10}.



Sections 2 and 3 of this paper are devoted to a detailed analysis of bounded positive solutions to the fully fractional heat equation \eqref{1.1}. The main focus lies on establishing monotonicity, symmetry, existence, and non-existence results under various structural assumptions on the functions \( a(x) \) and \( f(u) \).

In Section 2, under the assumptions that \( a(x) \) is strictly increasing in the \( x_1 \)-direction and \( f(u) \) is locally Lipschitz, non-decreasing on \( (0, +\infty) \), and satisfies \( f(0) = 0 \), \( f'(0) = 0 \), we prove that any uniformly continuous bounded positive solution \( u(x,t) \) must be strictly increasing in the \( x_1 \)-direction for each fixed time \( t \). This is achieved via the method of moving planes, which involves studying the reflected function \( w_\lambda(x,t) = u(x^\lambda, t) - u(x,t) \) and deriving a differential inequality for \( w_\lambda \) under the action of the fractional operator \( (\partial_t - \Delta)^s \). Furthermore, if \( a(x) \) exhibits certain symmetries (such as being even in \( x_1 \) or radially symmetric about the origin), we show that the solution \( u \) inherits the same symmetry properties.

Section 3 addresses the existence and non-existence of bounded positive solutions to the fully fractional heat equation \cite{1.1}. We first demonstrate that if \(a(x) \to \infty\) as \(x_1 \to \infty\), then no bounded positive solution exists. This is proved by constructing a suitable auxiliary function and applying a maximum principle argument to derive a contradiction. Furthermore, when \(a(x) \sim |x_1|^\alpha\) and \(f(u) \sim u^r\) as \(|x|, u \to \infty\), a critical exponent is identified:
\[
r^* = \min\left\{1 + \frac{4s^2}{n}, \frac{n + 2 + 2s - 2\alpha}{n + 2 - 2s}\right\}.
\]
In the subcritical case \(1 < r < r^*\), we establish the existence of a nontrivial bounded decaying solution via fractional heat semigroup theory and the Banach fixed point theorem. Compactness is verified using a priori estimates and the singular Gronwall inequality. In the critical and supercritical cases \(r \geq r^*\), we prove that no such solution exists, by showing that the necessary scaling identity derived from the energy functional cannot hold, and further demonstrate that finite-time blow-up occurs for certain initial data through a concavity argument and test function method.


These results not only extend classical theories for local operators to the nonlocal setting but also introduce new techniques for handling the fully fractional heat operator, which may be applicable to other nonlocal problems.

\section{Monotonicity and Symmetry of Bounded Positive Solutions}

\quad \quad This section will prove the monotonicity and symmetry of bounded positive solutions to equation \eqref{1.1} under certain conditions.

\textbf{Theorem 2.1.} Let \(u(x,t) \in C^{2s+\epsilon,s+\epsilon}_{x,t,loc}(\mathbb{R}^n \times \mathbb{R}) \cap L(\mathbb{R}^n \times \mathbb{R})\) be a bounded classical positive solution of equation \eqref{1.1}, where \(0 < s < 1\). If \(u\) is uniformly continuous, and the following assumptions hold for \(a(x)\) and \(f(u)\):\newline
(H1) \(f(\cdot)\) is locally Lipschitz continuous, non-decreasing on \((0,+\infty)\), \(f'(\cdot)\) is continuous near zero, and
\begin{equation}
f(0) = 0, f'(0) = 0, f(\cdot) > 0 \text{ on } (0,+\infty). \label{2.1}\end{equation}
(H2)
\begin{equation}
a(x_1,x') \in C(\mathbb{R}^n) \text{ is strictly increasing in } x_1, \forall x' \in \mathbb{R}^{n-1}.\label{2.2}\end{equation}
and \(a(x_1,\cdot) \in L^\infty(\mathbb{R}^{n-1})\), furthermore,
\[\limsup_{x_1 \to -\infty} a(x_1,x')|x_1|^{2s} \leq 0, \forall x' \in \mathbb{R}^{n-1}.\]
Then for any \(t \in \mathbb{R}\), \(u(\cdot,t)\) is strictly increasing in the \(x_1\) direction.

\textbf{Proof.} The proof follows a similar line of reasoning as presented by \cite{28, 27} and is included here for  completeness.  We will apply the method of moving planes for the proof. First, we denote the following notations:
\begin{itemize}
    \item \(T_\lambda := \{x \in \mathbb{R}^n | x_1 = \lambda\}\) is the moving plane;
    \item \(\Sigma_\lambda := \{x \in \mathbb{R}^n | x_1 < \lambda\}\) is the region to the left of the moving plane;
    \item \(x^\lambda := (2\lambda - x_1, x_2, \cdots, x_n)\) is the reflection point with respect to the hyperplane \(T_\lambda\);
    \item \( u_{\lambda}(x,t) := u(x^{\lambda},t) \);
    \item \( w_{\lambda}(x,t) := u_{\lambda}(x,t) - u(x,t) \).
\end{itemize}
According to assumptions (H1) and (H2) and the boundedness of \(u,\) the following differential inequality is obtained: 
\begin{align}
(\partial_t - \Delta)^s w_{\lambda}(x,t)  &= a(x^\lambda)f(u_{\lambda}(x,t)) - a(x)f(u(x,t)) \nonumber \\
&= [a(x^\lambda) - a(x)] f(u_{\lambda}(x,t)) + a(x)[f(u_{\lambda}(x,t)) - f(u(x,t))] \nonumber \\
&\geq a(x)M_{\lambda}(x,t)w_{\lambda}(x,t), \label{2.3}
\end{align}
where
\[ M_{\lambda}(x,t) = \frac{f(u_{\lambda}(x,t)) - f(u(x,t))}{u_{\lambda}(x,t) - u(x,t)} \]
is bounded and non-negative. By the continuity and monotonicity of \(f,\) there exists a constant \(C>0,\) such that
\[ 0 \leq M_{\lambda}(x,t) \leq C. \]

\textbf{Step 1.} Prove that when \(\lambda\) is sufficiently negative, the following inequality holds:
\begin{equation}
w_{\lambda}(x,t) \geq 0, \forall (x,t) \in \Sigma_{\lambda} \times \mathbb{R}. \label{2.4} 
\end{equation}
We will prove it by contradiction. Assume that \eqref{2.4} does not hold. Since \(\lambda\) is sufficiently negative, there exists an undetermined constant \(R > 0\) and \(\lambda < -R\). Since \(u\) is bounded, let
\begin{equation}
A := -\inf_{(x,t) \in \Sigma_{\lambda} \times \mathbb{R}} w_{\lambda}(x,t) > 0. \label{2.5} 
\end{equation}
Since \(\Sigma_{\lambda}\) and \(\mathbb{R}\) are unbounded, there exists an approximate minimizing sequence \(((x^k, t_k)) \subset \Sigma_{\lambda} \times \mathbb{R}\) and a sequence \(\{\epsilon_k\}\) monotonically decreasing to 0, such that
\begin{equation}
w_{\lambda}(x^k, t_k) = -A + \epsilon_k < 0. \label{2.6} 
\end{equation}
To construct a sequence of functions achieving the minimum, perturb \(w_{\lambda}\) near \((x^k, t_k)\) as follows:
\begin{equation}
v_k(x,t) = w_{\lambda}(x,t) - \epsilon_k \eta_k(x,t), (x,t) \in U_{\delta}(x^k, t_k), \text{ for some } \delta > 0, \label{2.7} \end{equation}
where
\begin{equation}
\eta_k(x,t) = \eta \left( \frac{x - x^k}{r_k}, \frac{t - t_k}{r_k^2} \right), \label{2.8} \end{equation}
\[ r_k = \frac{1}{2}dist(x^k, T_\lambda) > 0,\]
\( \eta \in C_0^\infty(\mathbb{R}^n \times \mathbb{R}) \) is a cutoff function satisfying:
\[\begin{cases}
0 \leq \eta \leq 1, & (x,t) \in \mathbb{R}^n \times \mathbb{R}, \\
\eta = 1, & (x,t) \in B_{\frac{1}{2}}(0) \times \left[ -\frac{1}{2}, \frac{1}{2} \right], \\
\eta = 0, & (x,t) \in (\mathbb{R}^n \times (-\infty, 1]) \setminus (B_1(0) \times [-1, 1]).
\end{cases}\]
Define the parabolic cylinder centered at \((x^k, t_k)\) as
\[ Q_k(x^k, t_k) := B_{r_k}(x^k) \times \left[ t_k - r_k^2, t_k + r_k^2 \right] \subset \Sigma_\lambda \times \mathbb{R}. \]
From \eqref{2.5}-\eqref{2.7} we obtain
\[ v_k(x^k, t_k) = -A, \]
\[ v_k(x, t) = w_\lambda(x, t) \geq -A, (x, t) \in (\Sigma_\lambda \times \mathbb{R}) \setminus Q_k(x^k, t_k), \]
\[ v_k(x, t) = -\epsilon_k \eta_k(x, t) > -A, (x, t) \in T_\lambda \times \mathbb{R}. \]
This shows that each \(v_k\) must achieve its minimum at some point \((x^k, t_k)\) within the parabolic cylinder \(\overline{Q_k(x^k, t_k)} \subset \Sigma_\lambda \times \mathbb{R}\), and this minimum does not exceed \(-A\). From \eqref{2.5} and \eqref{2.7}, \(\exists (\bar{x}^k, \bar{t}_k) \in \overline{Q_k(x^k, t_k)}\), such that
\begin{equation} 
-A - \epsilon_k \leq v_k(\bar{x}^k, \bar{t}_k) = \inf_{\Sigma_\lambda \times \mathbb{R}} v_k(x, t) \leq -A. \label{2.9} \end{equation}
From this, it follows that
\begin{equation} -A \leq w_\lambda(\bar{x}^k, \bar{t}_k) \leq -A + \epsilon_k < 0. \label{2.10} \end{equation}
Furthermore, from the definition of the operator \((\partial_t - \Delta)^s\), the antisymmetry of \(w_\lambda\) with respect to \(x_1\), and for any \(y \in \Sigma_\lambda\),we have \(|x^k - y^\lambda| > |x^k - y|\), combined with \eqref{2.9}, we obtain
\begin{align}
(\partial_t - \Delta)^s v_k(\bar{x}^k, \bar{t}_k) &= C_{n,s} \int_{-\infty}^{\bar{t}_k} \int_{\mathbb{R}^n} \frac{v_k(\bar{x}^k, \bar{t}_k) - v_k(y, \tau)}{(\bar{t}_k - \tau)^{\frac{n}{2}+1+s}} e^{-\frac{|x^k-y|^2}{4(\bar{t}_k-\tau)}} dy d\tau \nonumber \\
&= C_{n,s} \int_{-\infty}^{\bar{t}_k} \int_{\Sigma_\lambda} \frac{v_k(\bar{x}^k, \bar{t}_k) - v_k(y, \tau)}{(\bar{t}_k - \tau)^{\frac{n}{2}+1+s}} e^{-\frac{|x^k-y|^2}{4(\bar{t}_k-\tau)}} dy d\tau \nonumber \\
&+ C_{n,s} \int_{-\infty}^{\bar{t}_k} \int_{\Sigma_\lambda} \frac{v_k(\bar{x}^k, \bar{t}_k) - v_k(y^\lambda, \tau)}{(\bar{t}_k - \tau)^{\frac{n}{2}+1+s}} e^{-\frac{|x^k-y|^2}{4(\bar{t}_k-\tau)}} dy d\tau \nonumber \\
&\leq C_{n,s} \int_{-\infty}^{\bar{t}_k} \int_{\Sigma_\lambda} \frac{2v_k(\bar{x}^k, \bar{t}_k) - v_k(y, \tau) - v_k(y^\lambda, \tau)}{(\bar{t}_k - \tau)^{\frac{n}{2}+1+s}} e^{-\frac{|x^k-y|^2}{4(\bar{t}_k-\tau)}} dy d\tau \nonumber \\
&\leq C_{n,s} 2\left(v_k(\bar{x}^k, \bar{t}_k) + \epsilon_k\right) \int_{-\infty}^{\bar{t}_k} \int_{\Sigma_\lambda} \frac{1}{(\bar{t}_k - \tau)^{\frac{n}{2}+1+s}} e^{-\frac{|x^k-y|^2}{4(\bar{t}_k-\tau)}} dy d\tau \nonumber \\
&\leq \frac{C_1(-A + \epsilon_k)}{r_k^{2s}}. \label{2.11}
\end{align}
Then, combining \eqref{2.7} and \eqref{2.11}, and since \((\partial_t - \Delta)^s \eta_k(\bar{x}^k, \bar{t}_k) \leq \frac{C_2}{r_k^{2s}}\), we have
\begin{align}
(\partial_t - \Delta)^s w_\lambda(\bar{x}^k, \bar{t}_k) &= (\partial_t - \Delta)^s v_k(\bar{x}^k, \bar{t}_k) + \epsilon_k (\partial_t - \Delta)^s \eta_k(\bar{x}^k, \bar{t}_k) \nonumber \\
&\leq \frac{C_1(-A + \epsilon_k)}{r_k^{2s}} + \frac{C_2 \epsilon_k}{r_k^{2s}} \nonumber \\
&\leq \frac{C(-A + \epsilon_k)}{r_k^{2s}}. \nonumber
\end{align}
Combining \(\lambda \leq 0\), inequality \eqref{2.10} and the differential inequality \eqref{2.3} yields
\begin{equation}
\frac{C(-A + \epsilon_k)}{r_k^{2s}} \geq a(x^k) M_\lambda(\bar{x}^k, \bar{t}_k) w_\lambda(\bar{x}^k, \bar{t}_k).\label{2.12}
\end{equation}
It is known that \(\frac{C(-A+\epsilon_k)}{r_k^{2s}} < 0\), and \(M_\lambda(\bar{x}^k, \bar{t}_k) \in [0,C]\), \(w_\lambda(\bar{x}^k, \bar{t}_k) < 0\). Thus, the sign of \(a(x^k)\) needs to be discussed.

If \(a(x^k) \leq 0\), then the right-hand side of \eqref{2.12} is greater than or equal to 0, contradicting \(\frac{C(-A+\epsilon_k)}{r_k^{2s}} < 0\); if \(a(x^k) > 0\), then from \eqref{2.5}, \eqref{2.12} and \(0 < M_\lambda \leq 1\), we obtain
\[\frac{C(-A + \epsilon_k)}{-A} \leq a(x^k) r_k^{2s} \leq a(x^k) |x^k|^{2s}.\]
Then for sufficiently large \(k\), we have
\[a(x^k) |x^k|^{2s} \geq \frac{C}{2},\]
which contradicts assumption (H2). In summary, we obtain that \eqref{2.4} holds.

\textbf{Step 2.} Inequality \eqref{2.4} provides the starting point for moving the plane. Now move the plane \(T_\lambda\) to the right along the x direction and maintain inequality \eqref{2.4} during the process until its limiting position \(T_{\lambda_0}\), where \(\lambda_0\) is defined as:
\[ \lambda_0 := \sup \left\{ \lambda : w_{\mu}(x,t) \geq 0, (x,t) \in \Sigma_{\mu} \times \mathbb{R}, \forall \mu \leq \lambda \right\}. \]
Now prove
\[ \lambda_0 = +\infty. \]
If not, then \(\lambda_0 < +\infty\), so by its definition, there exists a sequence \(\{\lambda_k\}\) monotonically decreasing to \(\lambda_0\) and some positive sequence \(\{m_k\}\), such that
\[ \inf_{(x,t) \in \Sigma_{\mu} \times \mathbb{R}} w_{\lambda_k}(x,t) := -m_k < 0.  \]

We first prove
\begin{equation}
\lim_{k \to \infty} m_k = 0. \label{2.13} 
\end{equation}
If \eqref{2.13} does not hold, then there exists a subsequence of \(\{m_k\}\) (still denoted by \(\{m_k\}\)), such that for some \(M > 0\), \(m_k > M\) (i.e., \(-m_k < -M\)). Thus, there exists a sequence \(\{(y^k, s_k)\} \subset \Sigma_{\lambda_k} \times \mathbb{R}\), such that
\begin{equation}
w_{\lambda_k}(y^k, s_k) \leq -M < 0. \label{2.14} 
\end{equation}
Now, consider cases based on the position of \(y^k\) to derive a contradiction with \eqref{2.14}.

If \(y^k \in \Sigma_{\lambda_k} \setminus \Sigma_{\lambda_0}\), then since \(\lambda_k \to \lambda_0\), when \(k \to \infty\),
\[ |y^k - (y^k)^{\lambda_k}| = 2|\lambda_k - y^k_1| \to 0. \]
Therefore, by the uniform continuity of u, when \(k \to \infty\),
\[ w_{\lambda_k}(y^k, s_k) = u\left((y^k)^{\lambda_k}, s_k\right) - u\left(y^k, s_k\right) \to 0. \]
This contradicts \eqref{2.14}.

If \( y^k \in \Sigma_{\lambda_0} \), then combining \(\lambda_k \to \lambda_0\), the uniform continuity of \(u\), and the definition of \(\lambda_0\), we get when \( k \to \infty \),
\[w_{\lambda_k}(y^k, s_k) = u((y^k)^{\lambda_k}, s_k) - u((y^k)^{\lambda_0}, s_k) + w_{\lambda_0}(y^k, s_k)\]
\[\geq u((y^k)^{\lambda_k}, s_k) - u((y^k)^{\lambda_0}, s_k) \to 0.\]
This also contradicts \eqref{2.14}. Therefore \eqref{2.13} holds.

From \eqref{2.13}, there exists a sequence \((x^k, t_k) \subset \Sigma_{\lambda_k} \times \mathbb{R}\), such that
\[w_{\lambda_k}(x^k, t_k) = -m_k + m_k^2 < 0.\]
Perturb \( w_{\lambda_k} \) near \((x^k, t_k)\) as follows
\[v_k(x, t) = w_{\lambda_k}(x, t) - m_k^2 \eta_k(x, t), (x, t) \in \mathbb{R}^n \times \mathbb{R},\]
where \(\eta_k\) satisfies \eqref{2.8}, \(r_k = \frac{1}{2} dist(x^k, T_{\lambda_k})\). Define
\[P_k(x^k, t_k) := B_{r_k}(x^k) \times [t_k - r_k^2, t_k + r_k^2] \subset \Sigma_{\lambda_k} \times \mathbb{R},\]
then each \( v_k \) achieves its minimum within \( P_k(x^k, t_k) \subset \Sigma_{\lambda_k} \times \mathbb{R} \), and it does not exceed \(-m_k\), i.e.:
\[\exists (\bar{x}^k, \bar{t}_k) \in P_k(x^k, t_k), \text{ such that } -m_k - m_k^2 \leq v_k(\bar{x}^k, \bar{t}_k) = \inf_{\Sigma_{\lambda_k} \times \mathbb{R}} v_k(x, t) \leq -m_k.\]
That is
\begin{equation}
-m_k \leq w_{\lambda_k}(\bar{x}^k, \bar{t}_k) \leq -m_k + m_k^2 < 0. \label{2.15}
\end{equation}
Through a process similar to Step 1, we obtain:
\begin{equation}
\frac{C(-m_k + m_k^2)}{r_k^{2s}} \geq a(x^k)M_{\lambda_k}(\bar{x}^k, \bar{t}_k)w_{\lambda_k}(\bar{x}^k, \bar{t}_k). \label{2.16}
\end{equation}

If \(a(x^k) \leq 0,\) then from \eqref{2.15}, we derive that the righ-hand side of \eqref{2.16} is equal to or greater than 0, contradicting the fact that \(\frac{C(-m_k+m_k^2)}{r_k^{2s}}<0.\) If \(a(x^k) > 0,\) we further obtain
\begin{equation}
C(1 - m_k) \leq a(x^k)r_k^{2s}M_{\lambda_k}(\bar{x}^k, \bar{t}_k) \label{2.17}
\end{equation}
from \eqref{2.15}. 
When \(k\) is large enough, there exists a positive number \(\epsilon,\) such that \(0<m_k\leq\epsilon,\) thus \(1-m_k \geq 1-\epsilon.\) That is to say, the left-hand side of \eqref{2.17} is bounded from below, thus the right-hand side of \eqref{2.17} is bounded from below. The boundedness of \(M_{\lambda_k}\) and the assumption (H2) guarantee that there exists a constant \(c>0,\) such that 

\begin{equation}
a(x^k), r_k^{2s}, M_{\lambda_k}(\bar{x}^k, \bar{t}_k) \geq c > 0. \label{2.18}
\end{equation}
Then from the definition of \(M_{\lambda_k},\) we have \(M_{\lambda_k}=f'\left(\xi_{\lambda_k}(\bar{x}^k, \bar{t}_k) \right),\) where \(\xi_{\lambda_k}(\bar{x}^k, \bar{t}_k)\) lies between \(u_{\lambda_k}(\bar{x}^k, \bar{t}_k)\) and \(u(\bar{x}^k, \bar{t}_k).\) Combining with \eqref{2.1} in (H1), we know that \(u(\bar{x}^k, \bar{t}_k)\) is bounded away from zero, that is
\begin{equation}
u(\bar{x}^k, \bar{t}_k) \geq C > 0, f\left(u(\bar{x}^k, \bar{t}_k)\right) \geq C > 0. \label{2.19}
\end{equation}
From the differential inequality \eqref{2.3}, \eqref{2.16} can be rewritten as:
\begin{align}
\frac{C(-m_k + m_k^2)}{r_k^{2s}} &\geq \left[ a\left((x^k)^{\lambda_k}\right) - a(x^k)\right]f\left(u^{\lambda_k}(\bar{x}^k, \bar{t}_k)\right) \nonumber \\
&+ a(x^k)M_{\lambda_k}(\bar{x}^k, \bar{t}_k)w_{\lambda_k}(\bar{x}^k, \bar{t}_k). \label{2.20}
\end{align}
Note that \(\left| x^k_1 - \lambda_k\right| \sim r_k=\frac{1}{2}\text{dist}\{x^k,T_{\lambda_k}\}\), \(x^k\) satisfies \(w_{\lambda_k}(x^k,t_k)=-m_k+m_k^2\), \(m_k \to 0\) as \(k \to \infty\), we come to the conclusion that \(r_k \to 0\), thus from \eqref{2.18} and \eqref{2.2} in (H2), we get
\begin{equation}
a\left((x^k)^{\lambda_k}\right) - a(x^k) \geq C > 0. \label{2.21}\end{equation}
Additionally, since \(w_\lambda(\bar{x}^k, \bar{t}_k) \to 0\) as \(k \to \infty\), from \eqref{2.19} we have
\begin{equation}
f\left(u^{\lambda_k}(\bar{x}^k, \bar{t}_k)\right) \geq C > 0. \label{2.22}\end{equation}
for sufficiently large \(k\). From \eqref{2.18}, \eqref{2.21} and \eqref{2.22}, the right-hand side of \eqref{2.20} is greater than 0; from \eqref{2.15}, the left-hand side of \eqref{2.20} is less than 0, leading to a contradiction. Hence, \(\lambda_0 = +\infty\).

\textbf{Step 3.} In the previous steps, we have proven
\[w_\lambda(x,t) \geq 0, \forall (x,t) \in \Sigma_\lambda \times \mathbb{R}, \forall \lambda \in \mathbb{R}.\]
Now prove that the strict inequality holds, i.e.,
\begin{equation}
w_\lambda(x,t) > 0, \forall (x,t) \in \Sigma_\lambda \times \mathbb{R}, \forall \lambda \in \mathbb{R}. \label{2.23}\end{equation}

If \eqref{2.23} does not hold, then for some fixed \(\lambda\), there exists a point \((x^0, t_0) \in \Sigma_\lambda \times \mathbb{R}\), such that
\begin{equation}
w_\lambda(x^0, t_0) = \min_{\Sigma_\lambda \times \mathbb{R}} w_\lambda(x,t) = 0. \label{2.24}
\end{equation}
From the differential inequality \eqref{2.3}, we get
\begin{equation}
(\partial_t - \Delta)^s w_\lambda(x^0, t_0) \geq 0. \label{2.25}
\end{equation}
On the other hand, at the minimum point \((x^0, t_0) \in \Sigma_\lambda \times \mathbb{R}\), it must hold that

\begin{align}
(\partial_t - \Delta)^s w_\lambda(x^0, t_0) &= C_{n,s} \int_{-\infty}^{t_0} \int_{\mathbb{R}^n} \frac{w_\lambda(x^0, t_0) - w_\lambda(y, \tau)}{(t_0 - \tau)^{\frac{n}{2}+1+s}} e^{-\frac{|x^0-y|^2}{4(t_0-\tau)}} dy d\tau \nonumber \\
&= C_{n,s} \int_{-\infty}^{t_0} \int_{\Sigma_\lambda} \frac{w_\lambda(x^0, t_0) - w_\lambda(y, \tau)}{(t_0 - \tau)^{\frac{n}{2}+1+s}} e^{-\frac{|x^0-y|^2}{4(t_0-\tau)}} dy d\tau \nonumber \\
&+ C_{n,s} \int_{-\infty}^{t_0} \int_{\Sigma_\lambda} \frac{w_\lambda(x^0, t_0) - w_\lambda(y^\lambda, \tau)}{(t_0 - \tau)^{\frac{n}{2}+1+s}} e^{-\frac{|x^0-y|^2}{4(t_0-\tau)}} dy d\tau \nonumber \\
&= C_{n,s} \int_{-\infty}^{t_0} \int_{\Sigma_\lambda} \frac{w_\lambda(x^0, t_0) - w_\lambda(y, \tau)}{(t_0 - \tau)^{\frac{n}{2}+1+s}} e^{-\frac{|x^0-y|^2}{4(t_0-\tau)}} dy d\tau \nonumber \\
&+ 2C_{n,s} w_\lambda(x^0, t_0) \int_{-\infty}^{t_0} \int_{\Sigma_\lambda} \frac{1}{(t_0 - \tau)^{\frac{n}{2}+1+s}} e^{-\frac{|x^0-y|^2}{4(t_0-\tau)}} dy d\tau \nonumber \\
&= C_{n,s} \int_{-\infty}^{t_0} \int_{\Sigma_\lambda} \frac{w_\lambda(x^0, t_0) - w_\lambda(y, \tau)}{(t_0 - \tau)^{\frac{n}{2}+1+s}} e^{-\frac{|x^0-y|^2}{4(t_0-\tau)}} dy d\tau \leq 0. \label{2.26}
\end{align}
From \eqref{2.25} and \eqref{2.26}, it follows that
\[C_{n,s} \int_{-\infty}^{t_0} \int_{\Sigma_\lambda} \frac{w_\lambda(x^0, t_0) - w_\lambda(y, \tau)}{(t_0 - \tau)^{\frac{n}{2}+1+s}} e^{-\frac{|x^0-y|^2}{4(t_0-\tau)}} dy d\tau = 0.\]
From \eqref{2.24}, \(w_\lambda(x^0, t_0) - w_\lambda(y, \tau) \leq 0\), so
\[w_\lambda(x^0, t_0) - w_\lambda(x, t) \equiv 0, \forall (x, t) \in \Sigma_\lambda \times (-\infty, t_0]\]
(Otherwise, the integral would not be 0, which leads to a contradiction). Combining again with \eqref{2.24}, we come to the following.
\[w_\lambda(x, t) \equiv 0, \forall (x, t) \in \Sigma_\lambda \times (-\infty, t_0].\]
Thus
\[u_\lambda(x, t) \equiv u(x, t), \forall (x, t) \in \mathbb{R}^n \times (-\infty, t_0].\]
This contradicts the second equality in \eqref{2.3}
\[(\partial_t - \Delta)^s w_\lambda(x, t) = [a(x^\lambda) - a(x)] f(u_\lambda(x, t)) + a(x) [f(u_\lambda(x, t)) - f(u(x, t))],\]
as the assumptions (H1) and (H2) on \(f(u)\) and \(a(x)\), and the positivity of \(u\), imply that the right-hand side of the above equation is greater than 0, while the left-hand side equals 0. Therefore, \eqref{2.23} holds, and \(u(x, t)\) is strictly increasing in the \(x_1\) direction.

\textbf{Remark 2.1.} In particular, the above conclusions hold when \( a(x) = x_1^k \), where \( k \) is a rational number with both numerator and denominator being positive odd integers; and \( f(u) = u^p \), where \( 1 < p < +\infty \).

\textbf{Theorem 2.2.} Let \( u(x,t) \in C_{x,t,loc}^{2s+\epsilon,s+\epsilon} (\mathbb{R}^n \times \mathbb{R}) \cap L(\mathbb{R}^n \times \mathbb{R}) \) be a bounded classical positive solution of equation \eqref{1.1}, \( 0 < s < 1 \). Suppose \( u \) is uniformly continuous, \( f(u) \) satisfies assumption (H1) in Theorem 2.1, and \( a(x) \) satisfies the following assumption:\newline
(H2')
\[a(x_1,x') \in C(\mathbb{R}^n) \text{ is strictly increasing in } x_1 \text{ when } x_1 > 0, \forall x' \in \mathbb{R}^{n-1}\]
and \( a(x_1, \cdot) \in L^\infty(\mathbb{R}^{n-1}) \), furthermore, \( a(x) \) is an even function in \( x_1 \) or radially symmetric about the origin.
Then\newline
(i) For any \( t \in \mathbb{R}, u(x,t) \) has the same symmetry as \( a(x) \);\newline
(ii) It is strictly increasing in the \( x_1 \) direction when \( x_1 < 0 \).

\textbf{Proof.} 
(i) Assume \( u \) is a solution of equation \eqref{1.1}, i.e.,
\[(\partial_t - \Delta)^s u = a(x)f(u).\]
We will show that when \( a(x) \) is even in \( x_1 \), \( u_0(x,t) := u(x^0,t) := u(-x_1,x',t) \) is also a solution of equation \eqref{1.1}.

Define the reflection operator \( PU(x_1,x_2,\cdots,x_n,t) = U(-x_1,x_2,\cdots,x_n,t) \), then \( Pu(x,t) = u_0(x,t) \). Since \( a(x) \) is even in \( x_1 \), we have:
\[Pa(x) = a(-x_1,x') = a(x).\]
Thus, the right-hand side of the equation
\begin{equation}
P(a(x)f(u)) = a(x)f(u_0) = a(x)f(u(-x_1,x',t)) = a(x)f(Pu). \label{2.27}\end{equation}
Through Fourier transform\cite{25, 15, 26}, it can be shown that the operator \( (\partial_t - \Delta)^s \) commutes with the reflection operator \( P \), i.e.,
\begin{equation}
P((\partial_t - \Delta)^s u(x,t)) = (\partial_t - \Delta)^s (Pu(x,t)). \label{2.28}\end{equation}
This is because in Fourier space, reflection only changes the sign of \(\xi_1\), while the symbol of the operator \(\{i\tau + |\xi|^2\}\) only depends on \(|\xi|^2\). Combining \eqref{1.1}, \eqref{2.27} and \eqref{2.28} yields
\begin{equation}
(\partial_t - \Delta)^s (Pu(x,t)) = a(x)f(Ru).\label{2.29}\end{equation}

From the definition of the full fractional heat operator, we know that
\begin{align}
(\partial_t - \Delta)^s u(x^0, t) &= C_{n,s} \int_{-\infty}^t \int_{\mathbb{R}^n} \frac{u(x^0, t) - u(y, \tau)}{(t - \tau)^{\frac{n}{2}+1+s}} e^{-\frac{|x^0 - y|^2}{4(t - \tau)}} dy d\tau \nonumber \\
&= C_{n,s} \int_{-\infty}^t \int_{\mathbb{R}^n} \frac{u(x^0, t) - u(y^0, \tau)}{(t - \tau)^{\frac{n}{2}+1+s}} e^{-\frac{|x^0 - y^0|^2}{4(t - \tau)}} dy^0 d\tau \nonumber \\
&= C_{n,s} \int_{-\infty}^t \int_{\mathbb{R}^n} \frac{u(x, t) - u(y, \tau)}{(t - \tau)^{\frac{n}{2}+1+s}} e^{-\frac{|x^0 - y^0|^2}{4(t - \tau)}} dy d\tau \nonumber \\
&= (\partial_t - \Delta)^s u(x, t). \label{2.30}
\end{align}
where \(x^0 = (-x_1, x'), y^0 = (-y_1, y')\). The second equality follows from the arbitrariness of \(y\), and the third equality is obtained by the change of variable \(y = y^0\) (its Jacobian determinant is 1). Thus, \((\partial_t - \Delta)^s u_0(x, t) = (\partial_t - \Delta)^s u(x, t)\). Combining \eqref{2.29} and \eqref{2.30} yields \(Ru = u\), i.e., u is even in \(x_1\). When \(a(x)\) is radially symmetric about the origin, it can be similarly shown that u is radially symmetric about the origin. Therefore, the solution \(u(x, t)\) has the same symmetry as \(a(x)\).

(ii) The proof of this part is similar to the proof of part (ii) in Theorem 2.1. The main differences are replacing the assumption (H2) involved in the process with (H2'), and changing the definition of the supremum \(\lambda\) mentioned in Step 2 to
\[\lambda_0 := \sup \left\{ \lambda \leq 0 : w_\mu(x, t) \geq 0, (x, t) \in \Sigma_\mu \times \mathbb{R}, \forall \mu \leq \lambda \right\}\]
and proving its value is 0, where the sequence \(\{\lambda_k\}\) constructed in the contradiction argument satisfies \(\lambda_0 < \lambda_k < 0\) and \(\{\lambda_k\}\) monotonically decreases to \(\lambda_0\).

\section{Existence and Non-existence of Bounded Positive Solutions}
\quad \quad This section mainly proves the following theorem:

\textbf{Theorem 3.1.} Based on Theorem 2.1, further assume that for any \(x' \in \mathbb{R}^{n-1}\),
\begin{equation}
\lim_{x_1 \to \infty} a(x_1, x') = \infty, \label{3.1}\end{equation}
then equation \eqref{1.1} does not exist a bounded positive solution.

\textbf{Proof.} Let \(\lambda_1\) and \(\phi\) be the first eigenvalue and the corresponding eigenfunction of the problem
\[\begin{cases}
(-\Delta)^s \phi(x) = \lambda_1 \phi(x), & x \in B_1(0)\\
\phi(x) = 0, & x \in B_1^c(0)
\end{cases}\]
and assume
\[\max_{\mathbb{R}^n} \phi(x) = 1.\]
For any \(R \geq 1\), define
\[\phi_R(x) = \phi(x - Re_1),\]
where \(e_1\) is the unit vector in the \(x_1\) direction.

Let
\[v(x,t) = \phi_R(x) \eta(t),\]
where \(\eta(t) = t^\beta - 1\), \(0 < \beta = \frac{1}{2k+1} < s\), and \(k\) is some positive integer.

We prove the theorem by contradiction. Assume \(u\) is a bounded positive solution of \eqref{1.1}. Define
\[T = (M + 1)^{\frac{1}{\beta}}, M := \sup_{\mathbb{R}^n \times \mathbb{R}} u(x,t).\]
Thus \(M = T^\beta - 1\). Then for any \((x,t) \in B_1(Re_1) \times [1,T]\), we have
\begin{align}
(\partial_t - \Delta)^s v(x,t) &= C_{n,s} \int_{-\infty}^t \int_{\mathbb{R}^n} \frac{\phi_R(x) \eta(t) - \phi_R(y) \eta(\tau)}{(t - \tau)^{\frac{n}{2}+1+s}} e^{-\frac{|x-y|^2}{4(t-\tau)}} dy d\tau \nonumber \\
&= C_{n,s} \int_{-\infty}^t \int_{\mathbb{R}^n} \frac{[\phi_R(x) - \phi_R(y)] \eta(t) - \phi_R(y)[\eta(t) - \eta(\tau)]}{(t - \tau)^{\frac{n}{2}+1+s}} e^{-\frac{|x-y|^2}{4(t-\tau)}} dy d\tau \nonumber \\
&\leq \eta(t)(-\Delta)^s \phi_R(x) + \sup_{\mathbb{R}^n} \phi_R(x) \partial_t^s \eta(t) \nonumber \\
&\leq (\lambda_1 \eta(t) + C_s t^{\beta-s}) \sup_{\mathbb{R}^n} \phi_R(x) \nonumber \\
&\leq \lambda_1(T^\beta - 1) + C_s := C_T. \label{3.2}
\end{align}
On the other hand, from equation \eqref{1.1}, for any \((x,t) \in B_1(Re_1) \times [1,T]\), we have
\begin{align}
(\partial_t - \Delta)^s u(x,t) &\geq a(R-1,x') f(u(x,t)) \nonumber \\
&\geq a(R-1,x') f(m_T),\label{3.3}
\end{align}
where
\[ m_T := \inf_{B_1(Re_1) \times [1, T]} u(x, t) > 0. \]
 
Using the monotonicity in \( x_1 \) from Theorem 2.1. From condition \eqref{3.1}, we can choose sufficiently large \( R \), such that
\begin{equation}
a(R - 1, x')f(m_T) > C_T. \label{3.4} 
\end{equation}
Let
\[ w(x, t) = u(x, t) - v(x, t), \]
Then from \eqref{3.2}-\eqref{3.4}, \( w \) satisfies:
\[\begin{cases}
(\partial_t - \Delta)^s w(x, t) \geq 0, & (x, t) \in B_1(Re_1) \times (1, T], \\
w(x, t) > 0, & (x, t) \in (\mathbb{R}^n \times (-\infty, T]) \setminus (B_1(Re_1) \times (1, T]). 
\end{cases} \]
Thus, by the maximum principle,
\[ w(x, t) > 0, \forall (x, t) \in B_1(Re_1) \times (1, T]. \]
That is,
\[ u(x, t) > \phi_R(x)(t^\beta - 1), \forall (x, t) \in B_1(Re_1) \times (1, T], \]
Hence
\[ M > \max_{\mathbb{R}^n} \phi_R(x)(T^\beta - 1) = 1 \cdot M = M, \]
leading to a contradiction. Therefore, equation \eqref{1.1} does not have a bounded positive solution.

\textbf{Theorem 3.2.} For the full fractional heat operator equation
\begin{equation}\label{3.5}
    (\partial_t - \Delta)^s u(x,t) = a(x)f(u), \quad (x,t) \in \mathbb{R}^n \times (0,+\infty),
\end{equation}
where \(a(x)\) satisfies \(|a(x)| \leq C(1 + |x_1|^\alpha)\), \(f(u)\) satisfies \(|f(u)| \leq C(1 + |u|^\gamma)\) and is locally Lipschitz continuous, with initial data \(u_0(x) = u(x,0) \in L^2(\mathbb{R}^n)\), then:\\ 
(i)There exists \(T > 0\) and a unique function
\[u \in C\left([0,T]; L^2(\mathbb{R}^n)\right) \cap L^2(0,T; H^s(\mathbb{R}^n))\]
satisfying the integral equation
\begin{equation}\label{3.5.5}
u(t) = S(t)u_0(x) + \int_0^t K(t-\tau)[a(x)f(u(x,\tau))]d\tau,\end{equation}
where \(S(t) = e^{-(t-\Delta)^s}\) is the fractional heat semigroup and \(K(t)\) is the corresponding kernel function. If it is further assumed that \(1 < r < r^*\), then the local solution can be extended to a global solution:
\[u \in C([0, +\infty); L^2(\mathbb{R}^n)) \cap L_{\text{loc}}^2(0, +\infty; H^s(\mathbb{R}^n)).\]
Furthermore, for any \(\sigma < 2s^2\), we have
\[u \in C(0, +\infty; H^{s+\sigma}(\mathbb{R}^n)).\]
Here, the critical exponent \(r^*\) is defined as
\[r^*=\min\left\{\ 1+\frac{4s^2}{n},\frac{n+2+2s-2\alpha}{n+2-2s} \right\}.\]
(ii)(Blow-up Phenomena) Assume \(r \geq r^*\), and there exists a constant \(C > 0\) such that
\[a(x) \geq C|x_1|^{\alpha}, \quad f(u) \geq C|u|^{r}\]
hold for all \(x \in R^n, u \in R\). Then there exists initial data \(u_0 \in L^2(R^n)\) such that the solution to equation \eqref{3.5} blows up in finite time.

\textbf{Proof.} Define the space-time function space
\[X_T = \left\{u \in C\left([0,T]; L^2(\mathbb{R}^n)\right): \|u\|_{X_T} < \infty\right\},\]
where the norm is
\[\|u\|_{X_T} = \sup_{0 \leq t \leq T}\|u(x,t)\|_{L^2} + \left(\int_0^T \|u(x,t)\|_{H^s}^{2}dt\right)^{1/2}.\]
Define the solution operator \(\Phi: X_T \to X_T\) as
\[\Phi(u) = S(t)u_0 + \int_0^t K(t-\tau)[a(x)f(u(x,\tau))]d\tau,\]
then the problem reduces to finding a fixed point of \(\Phi: X_T \to X_T\). By the theory of the fractional heat semigroup\cite{30}, the kernel function \(K(t)\) satisfies the estimate:
\begin{equation}\label{3.6}
\|K(t)f\|_{L^2} \leq Ct^{s-1}\|f\|_{L^2}, \quad t > 0.
\end{equation}
Then for any \(u \in X_T\), we have
\begin{align}\label{3.7}
\|\Phi(u)(x,t)\|_{L^2} &\leq \|S(t)u_0\|_{L^2} + \int_0^t \|K(t-\tau)[a(x)f(u(x,\tau))]\|_{L^2}d\tau \nonumber \\
&\leq \|u_0\|_{L^2} + C \int_0^t (t-\tau)^{s-1}\|a(x)f(u(x,\tau))\|_{L^2}d\tau.
\end{align}
Using the growth conditions and weighted Sobolev embedding\cite{31}, there exists a constant \(C_\alpha > 0\) such that
\begin{equation}\label{3.8}
\|a(x)f(u)\|_{L^2} \leq C_\alpha \left(1 + \|u\|_{L^2}^\alpha\right).\end{equation}
Furthermore, by fractional Sobolev embedding and interpolation inequalities\cite{31,32}, when \( r < r^* \),
\begin{equation}\label{3.9}
\|u(t)\|_{L^2} \leq Ct^{-\beta} \sup_{0 \leq \tau \leq t} \|u(x, \tau)\|_{L^2}, \quad \beta = \frac{n}{4s} \left(1 - \frac{1}{r}\right).\end{equation}
Combining \eqref{3.6}-\eqref{3.9}, we obtain
\[\|\Phi(u)(x, t)\|_{L^2} \leq \|u_0(x)\|_{L^2} + Ct^s + Ct^{s-\beta} \|u\|_{X_T}^\sigma.\]
Similarly, it can be shown that \(\Phi(u) \in L^2(0, T; H^s)\). Choosing a sufficiently small \(T > 0\), it can be proven that \(\Phi\) is a contraction mapping on some closed ball in \(X_T\). The local existence and uniqueness of the solution then follow from the Banach fixed point theorem.

In order to extend local solutions to global solutions, we establish a uniform priori estimate. Define \(M(t) = \sup_{0 \leq \tau \leq t} \|u(x, \tau)\|_{L^2}\). From the integral equation \eqref{3.5.5}, we have:
\begin{equation}\label{3.10}
M(t) \leq \|u_0(x)\|_{L^2} + C \int_0^t (t - \tau)^{s-1} \left(1 + \tau^{-\beta} [M(\tau)]^r\right) d\tau.\end{equation}
\eqref{3.10} is a weakly singular Volterra inequality. By the singular Gronwall theory in \cite{33}, when \(s - \beta r > 0\) (i.e., \(r < r^*\)), there exists \(C(T) > 0\) such that \(M(t) \leq C(T)\) holds for all \(t \in [0, T]\). Furthermore, based on the uniform boundedness in \(L^\infty\), it can be further proven that
\[\int_0^T \|u(x, t)\|_{H^\infty}^2 dt \leq C(T).\]
This priori estimate ensures that the solution does not blow up in finite time, thus allowing the global solution to be obtained by extension from the local solution.

By the instantaneous regularization property of the fractional heat semigroup\cite{30}, for some \(\sigma > 0\), we have
\begin{equation}\label{3.11}
\|S(t)f\|_{H^\infty} \leq Ct^{-\frac{\sigma}{2s}} \|f\|_{L^2}.\end{equation}
From \eqref{3.5.5} and \eqref{3.11}, we have
\begin{equation}\label{3.12}
\|u(t)\|_{H^\infty} \leq Ct^{-\frac{\sigma}{2s}} \|u_0\|_{L^2} + C \int_0^t (t - \tau)^{s-1} \|u(x)f(u(x, \tau))\|_{L^2} d\tau.\end{equation}
When \(\sigma < 2s^2\), the singular integral in the above expression converges.

From the fact that \(u\) is a global solution, we have \(u \in L^\infty(0,T;L^2)\). Combining this with \eqref{3.12} yields \(u \in L^p(0,T;H^\sigma)\) for some \(p > 2\). Using this improved regularity, we find that \(a(x)f(u) \in L^q(0,T;L^2)\) for a larger \(q\). Repeating this process gradually enhances the regularity until the conclusion is reached.

(ii)Construct a family of test functions \(\phi_R(x,t)|_{R>0}\), where
\[\phi_R(x,t) = \psi\left(\frac{|x|}{R}\right)\eta\left(\frac{t}{R^{2s}}\right)\]
satisfying
\[\psi(s) =
\begin{cases}
1, & |s| \leq 1, \\
0, & |s| \geq 2;
\end{cases}\]
\[\eta(\tau) =
\begin{cases}
1, & |\tau| \leq 1, \\
0, & |\tau| \geq 2;
\end{cases}\]
with \(0 \leq \psi, \eta \leq 1\), uniformly compactly supported in space and time, and \(\psi'(s), |\eta'(\tau)|\) bounded.

Assume there exists a global solution \(u \in C\left([0,+\infty);L^2(R^n)\right)\). Multiply equation \eqref{3.5} by the test function \(\phi_R\) and integrate over the space-time domain:
\[\int_0^{+\infty} \int_{R^n} (\partial_t - \Delta)^s u \phi_R dxdt = \int_0^{+\infty} \int_{R^n} a(x)f(u)\phi_R dxdt.\]
Based on the self-adjointness of the full fractional heat operator\cite{1}, the left-hand side of the above equation can be transformed to obtain
\begin{equation}\label{3.14}
\int_0^{+\infty} \int_{R^n} u(\partial_t - \Delta)^s \phi_R dxdt = \int_0^{+\infty} \int_{R^n} a(x)f(u)\phi_R dxdt. \end{equation}
Perform the change of variables \(x = Ry, t = R^{2s}\tau\). Then the test function becomes
\[\phi_R(Ry, R^{2s}\tau) = \psi(|y|)\eta(\tau) = \phi(y,\tau),\]
The left-hand side of \eqref{3.14} becomes
\[\int_{0}^{\infty} \int_{\mathbb{R}^n} u(\partial_t - \Delta)^s \phi_R dxdt = R^{n+2s} \int_{0}^{\infty} \int_{\mathbb{R}^n} u(Ry, R^{2s}\tau)(\partial_t - \Delta)^s \phi(y, \tau)dyd\tau,\]
Using the lower bound condition on the nonlinear term, the right-hand side becomes
\begin{align}\label{3.15}
\int_{0}^{\infty} \int_{\mathbb{R}^n} a(x)f(u)\phi_R dxdt &\geq C \int_{0}^{\infty} \int_{\mathbb{R}^n} |x_1|^\alpha |u(x,t)|^r \phi_R (x,t)dxdt\nonumber \\
&= CR^{n+\alpha} \int_{0}^{\infty} \int_{\mathbb{R}^n} |y_1|^\alpha |u(Ry, R^{2s}\tau)|^r \phi(y, \tau)dyd\tau.\end{align}
Combining \eqref{3.14} and \eqref{3.15} yields
\begin{align}\label{3.16}
R^{n+2s} &\int_{0}^{\infty} \int_{\mathbb{R}^n} u(Ry, R^{2s}\tau)(\partial_t - \Delta)^s \phi(y, \tau)dyd\tau \nonumber \\ 
\geq CR^{n+\alpha} &\int_{0}^{\infty} \int_{\mathbb{R}^n} |y_1|^\alpha |u(Ry, R^{2s}\tau)|^r \phi(y, \tau)dyd\tau.\end{align}
Let \( U_R(y, \tau) = u(Ry, R^{2s}\tau) \), then \eqref{3.16} becomes
\begin{align}\label{3.17}
R^{n+2s} &\int_{0}^{\infty} \int_{\mathbb{R}^n} U_R(y, \tau)(\partial_t - \Delta)^s \phi(y, \tau)dyd\tau \nonumber \\
\geq CR^{n+\alpha} &\int_{0}^{\infty} \int_{\mathbb{R}^n} |y_1|^\alpha |U_R(y, \tau)|^r \phi(y, \tau)dyd\tau.\end{align}
Now assume the solution exists globally and does not blow up. We will show that a contradiction arises when \( R \to \infty \).

\textbf{Case 1.} When \( r > r^* \), rearranging \eqref{3.17} gives:
\[\int_{0}^{\infty} \int_{\mathbb{R}^n} |y_1|^\alpha |U_R(y, \tau)|^r \phi(y, \tau)dyd\tau \leq CR^{2s-\alpha-\gamma(r-1)},\]
where \(\gamma\) is the decay exponent of the solution. When \( r > r^* \), the exponent \( 2s - \alpha - \gamma(r-1) < 0 \). Letting \( R \to \infty \), the right-hand side tends to 0, while the left-hand side should remain positive if the solution is non-trivial, leading to a contradiction.

\textbf{Case 2.} When \( r = r^* \), define the auxiliary function
\[J(t) = \frac{1}{2} \int_{0}^{t} ||u(x, \tau)||_{L^2}^2 d\tau,\]
and compute its first and second derivatives:
\[J'(t) = \frac{1}{2} ||u(x, t)||_{L^2}^2,\]
\[J''(t) = \langle u(t), \partial_t u(t) \rangle.\]
Using the lower bound condition on the nonlinear term and the Hard-Littlewood-Sobolev inequality\cite{34}, it can be proven that there exists \( \varepsilon > 0 \) such that
\[J(t)J''(t) \geq (1 + \varepsilon) [J'(t)]^2.\]
The above inequality is equivalent to \(\frac{\text{d}}{\text{d}t} \left( \frac{J'(t)}{[J(t)]^{1+\varepsilon}} \right) > 0\), i.e., \(\frac{\text{d}^2}{\text{d}t^2} \left( [J(t)]^{-\varepsilon}\right) < 0\). This shows that \([J(t)]^{-\varepsilon}\) is a concave function.

Choose initial data \(u_0(x)\) and a constant \(C > 0\) depending on the parameters of the equation but independent of the solution, such that \(\|u_0(x)\|_{L^2}^2 > C\). Then by continuous dependence, there exists \(t_0 > 0\) such that \(J'(t_0) > 0\). For this \(t_0\), we have
\[\frac{J'(t)}{[J(t)]^{1+\varepsilon}} \geq \frac{J'(t_0)}{[J(t_0)]^{1+\varepsilon}}, \quad \text{i.e.,} \quad J'(t) \geq J'(t_0) \left( \frac{J(t)}{J(t_0)} \right)^{1+\varepsilon}\]
This causes \(J(t)\) to tend to infinity in finite time, meaning the solution blows up.

\bibliographystyle{plain}
\bibliography{ref}

Lu Haipeng

School of Mathematics and Statistics

Northwestern Polytechnical University

Xi'an, Shaanxi, 710129, P. R. China

luhp2023@mail.nwpu.edu.cn

Yu Mei

School of Mathematics and Statistics

Northwestern Polytechnical University

Xi'an, Shaanxi, 710129, P. R. China

yumei@nwpu.edu.cn

\end{document}